\documentclass{baart}
\usepackage{amssymb,amsmath,bbm,latexsym,amscd}
\usepackage[all]{xy}

\newtheorem{theorem}{Theorem}[section]
\newtheorem{proposition}[theorem]{Proposition}

\newtheorem{lemma}[theorem]{Lemma}
\newtheorem{corollary}[theorem]{Corollary}
\theoremstyle{definition}

\theoremstyle{remark}
\newtheorem{remark}[theorem]{Remark}

\numberwithin{equation}{section}

\begin{document}

 \title[Classification of Harish-Chandra modules]{Classification of Harish-Chandra Modules \\for Current Algebras}


\author{Michael Lau}
\address{D\'epartement de math\'ematiques et de statistique, Universit\'e Laval\\ Qu\'ebec, QC, Canada G1V 0A6}

\email{Michael.Lau@mat.ulaval.ca}
\thanks{Funding from the Natural Sciences and Engineering Research Council of Canada is gratefully acknowledged.}





\begin{abstract}
For any reductive Lie algebra $\mathfrak{g}$ and commutative, associative, unital algebra $S$, we give a complete classification of the simple weight modules of $\mathfrak{g}\otimes S $ with finite weight multiplicities.  In particular, any such module is parabolically induced from a simple admissible module for a Levi subalgebra.  Conversely, all modules obtained in this way have finite weight multiplicities.  These modules are isomorphic to tensor products of evaluation modules at distinct maximal ideals of $S$.  Our results also classify simple Harish-Chandra modules up to isomorphism for all central extensions of current algebras.

\end{abstract}

\maketitle





\renewcommand{\theequation}{\thesection.\arabic{equation}}
\newcommand{\K}{{k}}
\newcommand{\lra}{\longrightarrow}
\newcommand{\Q}{{\mathbb Q}}
\newcommand{\cN}{\mathcal{N}}
\newcommand{\cP}{\mathcal{P}}
\newcommand{\abG}{\left(\begin{array}{c}
a\\
b\end{array}\right)\ot \left(\begin{array}{c}
1\\
0\end{array}\right)}

\newcommand{\ghGbar}{\left(\begin{array}{c}
g\\
h\end{array}\right)\ot \left(\begin{array}{c}
0\\
1\end{array}\right)}

\newcommand{\Var}{\hbox{Var}\,}
\newcommand{\whp}{\widehat{\p}}
\newcommand{\wh}{\widehat}
\newcommand{\Curr}{\hbox{Curr}\,}
\newcommand{\Ctd}{\hbox{\text{\rm Ctd}}}
\newcommand{\ib}{{\overline{\iota}}}
\newcommand{\ot}{\otimes}
\newcommand{\ol}{\overline}
\newcommand{\oL}{\overline{\cL}}
\newcommand{\Z}{\mathbb{Z}}
\newcommand{\Cx}{\mathbb{C}}
\newcommand{\D}{\mathcal{D}}
\newcommand{\fF}{{\mathfrak F}}
\newcommand{\fG}{{\mathfrak G}}
\newcommand{\fg}{{\mathfrak g}}
\newcommand{\fh}{{\mathfrak h}}
\newcommand{\fE}{{\mathfrak E}}
\newcommand{\fl}{{\mathfrak l}}
\newcommand{\fp}{{\mathfrak p}}
\newcommand{\autfun}{{\bf Aut}}
\newcommand{\Rep}{\hbox{Rep}\,}
\newcommand{\Max}{\mathrm{Max}}
\newcommand{\p}{\partial}
\newcommand{\pAS}{\partial_{\A\ot\widehat{\cS}}}
\newcommand{\ga}{\alpha}
\newcommand{\gb}{\beta}
\newcommand{\gs}{\sigma}
\newcommand{\eps}{\epsilon}
\renewcommand{\proof}{{\bf Proof\ \ }}
\renewcommand{\qed}{\hfill $\Box$}
\newcommand{\bm}{{\bf m}}
\newcommand{\bmu}{{\mathbf \mu}}
\newcommand{\bsig}{{\mathbf \sigma}}
\newcommand{\bB}{{\bf B}}
\newcommand{\bG}{{\bf G}}
\newcommand{\bg}{{\bf g}}
\newcommand{\cL}{\mathcal{L}}

\newcommand{\gd}{\partial}
\newcommand{\gl}{\lambda}
\newcommand{\cU}{{\mathcal U}}
\newcommand{\cC}{{\mathcal C}}
\newcommand{\F}{{\mathcal F}}
\newcommand{\Vir}{{\mathcal Vir}}
\newcommand{\A}{{\mathcal A}}
\newcommand{\B}{{\mathcal B}}
\newcommand{\End}{\hbox{End}}
\newcommand{\pimap}[1]{\pi_{#1}}
\newcommand{\bo}[1]{\ol{\bf{#1}}}
\newcommand{\lb}[2]{\left\lbrack{#1}_\gl{#2}\right\rbrack}
\newcommand{\C}{{\mathbb C}}
\newcommand{\cG}{{\mathcal G}}
\newcommand{\R}{{\mathcal R}}
\newcommand{\cS}{{\mathcal S}}
\newcommand{\cT}{{\mathcal T}}
\newcommand{\even}{{\ol{0}}}
\newcommand{\odd}{{\ol{1}}}
\newcommand{\Aut}{{\rm Aut}}
\newcommand{\bAut}{{\bf Aut}}
\newcommand{\GL}{{\rm GL}}
\newcommand{\Hom}{{\rm Hom}}
\newcommand{\id}{{\rm id}}
\newcommand{\Id}{{\rm Id}}
\newcommand{\Span}{{\rm Span}}
\newcommand{\Ker}{{\rm Ker}}
\newcommand{\Mat}{{\rm Mat}}
\newcommand{\Gal}{{\mathcal Gal}}
\newcommand{\Loop}{{\mathcal L}}
\newcommand{\kalg}{{k\hbox{-alg}}}
\newcommand{\kalgder}{{k-\delta alg}}
\newcommand\limind{\mathop{\oalign{lim\cr
\hidewidth$\longrightarrow$\hidewidth\cr}}}
\newcommand\Spec{\text{\rm Spec}\,}
\newcommand\Out{\text{\rm Out}\,}
\newcommand\bOut{\text{\bf Out}\,}
\newcommand \bone{{\mathbf 1}}
\newcommand\Pic{\text{\rm Pic}}
\newcommand\os{\overset}
\newcommand\us{\underset}
\newcommand\et{\text{\rm \'et}}
\newcommand\ct{\text{\rm ct}}
\newcommand\fppf{\text{\rm fppf}}
\newcommand\q{\quad}
\newcommand{\SL}{{\mathbf{SL}}}
\newcommand{\bGL}{{\mathbf{GL}}}
\newcommand{\Zar}{{\rm Zar}}
\newcommand{\wtpsi}{\widetilde{\psi}}
\newcommand{\cI}{\mathcal{I}}
\newcommand{\cM}{\mathcal{M}}
\newcommand{\ev}{\mathrm{ev}}
\newcommand{\supp}{\hbox{supp}\,}
\def\mc{\mathcal}
\newcommand{\la}{\langle}
\newcommand{\ra}{\rangle}
\newcommand{\und}[1]{\underline{#1}}
\newcommand{\hgt}{\hbox{ht}\,}
\newcommand\blfootnote[1]{%
  \begingroup
  \renewcommand\thefootnote{}\footnote{#1}%
  \addtocounter{footnote}{-1}%
  \endgroup
}

\section{Introduction}\label{s1}  

The study of Harish-Chandra modules has a long and glorious history.  Two highlights were the description of composition factor multiplicities in the Kazhdan-Lusztig Conjectures, resolved by Beilinson-Bernstein and Brylinski-Kashiwara in 1981 \cite{BeBe81,BrKa81}, and the classification of simple Harish-Chandra modules for reductive Lie algebras, completed in 2000 by Mathieu, building on work of Fernando, Futorny, Benkart, Britten, and Lemire \cite{Ma00}.

In this paper, we classify the simple Harish-Chandra modules for generalised current algebras.  More precisely, let $\fg$ be a finite dimensional reductive Lie algebra over an algebraically closed field $k$ of characteristic $0$, and let $S$ be a commutative, associative, and unital $k$-algebra of finite type.  The associated {\em (generalised) current algebra} $\cG=\fg\ot_k S$ is the Lie algebra of all morphisms $\Spec\,S\rightarrow \fg$ of affine schemes.  Such algebras appear in various contexts of mathematics and physics, notably in affine Lie theory, sigma models, and string theory.  

A $\cG$-module $V$ is a {\em Harish-Chandra module} if it decomposes into finite dimensional common eigenspaces with respect to a Cartan subalgebra $\fh=\fh\ot 1$ of $\fg$.  Inspired by the classifications for reductive, affine, and Virasoro algebras \cite{Ma00,DiMaPe,DiGr09,FuTs01,Ma92,GLZ,Sa12}, we use locally nilpotent and injective actions of the elements of $\cG$ to prove a parabolic induction theorem, showing in Theorem \ref{thm1} that every simple Harish-Chandra module for a current algebra is parabolically induced from a weight module for a (generalised) Levi subalgebra $\cL$.  The absence of a $\Z$-grading in our very general context leads to different (and simpler) parabolic decompositions than in the affine and Virasoro cases.  We show that the relevant simple $\cL$-modules are {\em admissible}, meaning that their weight multiplicities are uniformly bounded.  Such representations were recently classified as tensor products of finitely many admissible evaluation modules \cite{curr}.  

The converse is not obvious.  The corresponding generalised Verma modules often have infinite dimensional weight spaces, so it is not a priori clear that their simple quotients, called {\em parabolically induced modules}, will always have finite weight multiplicities.  In Section 3, we prove finite dimensionality through a careful analysis of the action of raising operators on parabolically induced modules (Theorem \ref{thm2}).  Given any family of vectors in a weight space, the existence of a linear dependence is equivalent to a nontrivial solution of an infinite linear system determined by the vanishing of appropriate raising operators.  The challenge is to reduce this infinite family of equations to an underdetermined finite linear system, using properties of parabolic sets and exp-polynomial functions.  We obtain this reduction without using the Vandemonde determinant argument (which fails in our setting) introduced by Billig, Berman, and Zhao \cite{BeBi99,BiZh04} to settle the $\Z^n$-graded case.  

In Section 4, we present an alternative approach based on evaluation modules.  Using the classification of simple admissible modules for current algebras, we show that simple Harish-Chandra modules may be reinterpreted as tensor products of evaluation modules (Proposition \ref{star}).  Such $\cG$-modules are obtained by precomposing the action defining $\fg$-modules $V_i$ with reduction maps $\ev_i:\ \cG\rightarrow\fg$, $x\ot s\mapsto s(M_i)x$ at distinct maximal ideals $M_i$ of $S$.  This leads to a second proof of finite dimensionality of weight spaces, as well as an isomorphism criterion: such tensor products of evaluation modules will be isomorphic precisely when the corresponding sets of maximal ideals $M_i$ are equal and the corresponding $\fg$-modules $V_i$ are isomorphic.  Up to isomorphism, the classification thus reduces to the finite case, where $S$ is the base field $k$.  Using a trace argument, our results also extend to arbitrary central extensions of $\cG$.

\bigskip

\noindent
{\bf Acknowledgements.}  Much of this work was completed while the author was on sabbatical at the Institut Camille Jordan (Universit\'e de Lyon 1).  He thanks the institute for its warm hospitality during his visit.  The author also thanks the anonymous referee for very helpful comments and suggestions.

\section{Parabolic induction}\label{s2}

Throughout the paper, $\fh$ will denote a Cartan subalgebra of a finite dimensional reductive Lie algebra $\fg$ over an algebraically closed field $k$ of characteristic $0$.  All tensor products and algebras will be taken over the base field $k$ unless explicitly otherwise indicated.  We write $\mathbb{Z}$ and $\mathbb{N}$ for the set of integers and non-negative integers, respectively.  For each root $\ga$ in the root system $\Phi=\Phi(\fg,\fh)$, we fix a nonzero root vector $x_\ga$ and a Cartan element $h_\ga\in\fh$ so that $\{x_\ga,x_{-\ga},h_\ga\}$ forms an $\mathfrak{sl}_2$-triple, with $[x_\ga,x_{-\ga}]=h_\ga$ and $[h_\ga,x_\ga]=2x_\ga$.  For each commutative, associative, and unital $k$-algebra $S$ of finite type, the {\em (generalised) current algebra} is the Lie algebra $\cG=\fg\ot S$, with bracket $[x\ot r,y\ot s]=[x,y]\ot rs$, for all $x,y\in\fg$ and $r,s\in S$.  These algebras include the finite dimensional reductive Lie algebras $\fg$ and loop algebras $\fg\ot\Cx[t,t^{-1}]$ of affine Kac-Moody theory, as well as the $n$-point Lie algebras considered in classical and quantum field theory, multiloop algebras, Takiff algebras, and various Krichever-Novikov algebras.
  
We will consider {\em Harish-Chandra modules for $(\cG,\fh)$}, that is, $\cG$-modules $V$ which admit a direct sum decomposition into finite dimensional weight spaces with respect to the abelian subalgebra $\fh\ot 1\subseteq\cG$:
$$V=\bigoplus_{\lambda\in\fh^*}V_\lambda,$$
where $V_\lambda=\{v\in V\ :\ (h\ot 1)v=\lambda(h)v\hbox{\ for all\ }h\in\fh\}$ is finite dimensional for all $\lambda$.  
The set $\Pi(V)=\{\lambda\in\fh^*\ :\ V_\lambda\neq 0\}$ of {\em weights} is called the {\em support} of $V$.  
A Harish-Chandra module $V$ is {\em admissible} if the dimensions of its weight spaces $V_\lambda$ are uniformly bounded, that is, if there exists $N>0$ such that $\dim V_\lambda<N$ for all weights $\lambda$.

Note that if $\fg$ is reductive, then the simple Harish-Chandra modules of $\cG=\fg\ot S$ are tensor products of the simple Harish-Chandra modules of its components $\fg_i\ot S$, where $\fg=\oplus_{i=1}^n\fg_i$, and each $\fg_i$ is simple or abelian.  This can be proved by the argument in \cite[Proposition 3.4]{curr}, for instance.  Since the simple Harish-Chandra modules of an abelian Lie algebra are clearly one-dimensional, there is no loss of generality to assume that $\fg$ is simple.  {\em We will henceforth assume that $V$ is a simple Harish-Chandra module for $\cG=\fg\ot S$, where $\fg$ is finite dimensional and simple.}


\begin{lemma}\label{lemma1}
Let $\ga\in\Phi$ and $s\in S$.  Then either $x_\ga\ot s$ acts locally nilpotently everywhere on $V$, or else it acts injectively on $V$.
\end{lemma}

\noindent
\proof The proof is straightforward, using the fact that the set of vectors $v\in V$ on which $x_\ga\ot s$ acts locally nilpotently is a submodule of $V$.  This follows from the integrability of the Lie algebra $\fg$, viewed as an (adjoint) module over itself.\qed

\begin{proposition}\label{Phi^f-prop}
For any root $\ga\in\Phi$, the following conditions are equivalent:
\begin{enumerate}
\item[{\rm (i)}] For each $\lambda\in\Pi(V)$, the weight space $V_{\lambda+n\ga}$ is zero for all but finitely many $n>0$.
\item[{\rm (ii)}] There exists $\lambda\in\Pi(V)$ such that the weight space $V_{\lambda+n\ga}$ is zero for all but finitely many $n>0$.
\item[{\rm (iii)}] The element $x_\ga\ot s$ acts locally nilpotently on $V$ for all $s\in S$.
\item[{\rm (iv)}] The element $x_\ga\ot 1$ acts locally nilpotently on $V$.
\end{enumerate}
\end{proposition}

\noindent
\proof It is clear that (i) implies (ii), and (iii) implies (iv).  Moreover, if (ii) holds, then for all $s\in S$, we see that $x_\ga\ot s$ acts locally nilpotently on the nonzero space $V_\lambda$, and hence everywhere on $V$ by Lemma \ref{lemma1}.  Thus (ii) implies (iii).  

To see that (iv) implies (i), we recall an argument used by Fernando \cite{Fe90} in the finite dimensional context.  Namely, if $V_{\lambda+n\ga}$ is nonzero for infinitely many $n>0$, then there is a sequence of nonzero vectors $v_i\in V_{\lambda+n_i\ga}$ annihilated by $x_\ga\ot 1$, for positive integers $n_1<n_2<n_3<\cdots$.  Let $\mathfrak{s}_\ga=\hbox{Span}\{x_\ga\ot 1,x_{-\ga}\ot 1,h_\ga\ot 1\}$ be the copy of $\mathfrak{sl}_2(k)$ corresponding to the root $\ga$.  By taking $N\in\mathbb{Z}$ sufficiently large, we can guarantee that $(\lambda+N\ga)(h_\ga)$ is not a negative integer.  By $\mathfrak{sl}_2$-theory, the $\mathfrak{s}_\ga$-module $\mathcal{U}(\mathfrak{s}_\ga)v_i$ has a nonzero vector $w_i$ of weight $\lambda+N\ga$ for all $n_i>N$.  Each of the (infinitely many) modules $\mathcal{U}(\mathfrak{s}_\ga)v_i$ has a distinct highest weight, thus infinitely many of them have distinct central characters, and $\{w_i\ :\ n_i>N\}\subset V_{\lambda+N\ga}$ contains infinitely many linearly independent vectors.  This contradicts the finite dimensionality of the weight spaces of $V$.\qed

\bigskip

A root $\ga\in\Phi$ is called {\em locally finite} if one, and hence all, of the conditions of Proposition \ref{Phi^f-prop} holds.  We write $\Phi^f$ for the set of locally finite roots, and $\Phi^i=\Phi\setminus\Phi^f$ for its complement, the set of {\em injective} roots.  For any subset $T\subseteq \Phi$, let $-T$ be the set $\{-\ga\ :\ \ga\in T\}\subseteq \Phi$.  The set $P=\Phi^f\cup-\Phi^i$ will play a crucial role in what follows.

Consider the following subspaces of the Lie algebra $\fg$:
\begin{eqnarray*}
\mathfrak{l}&=&\fh\oplus\sum_{\ga\in P\cap-P}\fg_\ga\\
\mathfrak{n_\pm}&=&\sum_{\pm\ga\in P\setminus-P}\fg_\ga\\
\mathfrak{p}&=&\mathfrak{l}\oplus\mathfrak{n_+}.
\end{eqnarray*}
We will show that $\mathfrak{l}$, $\mathfrak{n_\pm}$, and $\mathfrak{p}$ are Lie subalgebras of $\fg$.  The critical step is the following lemma:
\begin{lemma}\label{lemma3}
\begin{enumerate}
\item[{\rm (i)}] Let $\ga,\beta\in\Phi^i$ and $\ga+\beta\in\Phi$.  Then $\ga+\beta\in\Phi^i$.
\item[{\rm (ii)}] Let $\ga,\beta\in\Phi^f$ and $\ga+\beta\in\Phi$.  Then $\ga+\beta\in\Phi^f$.
\item[{\rm (iii)}] Let $\ga\in\Phi^f\cap-\Phi^f$ and $\beta\in\Phi^i\cap-\Phi^i$.  Then $\ga+\beta\notin\Phi$.
\end{enumerate}
\end{lemma}

\noindent
\proof (i) By Lemma \ref{lemma1} and Proposition \ref{Phi^f-prop}, the elements $x_\ga\ot 1$ and $x_\beta\ot 1$ act injectively on the simple module $V$.  Thus $(x_\ga\ot 1)(x_\beta\ot 1)$ acts injectively, so $V_{\lambda+n(\ga+\beta)}$ is nonzero for all $\gl\in\Pi(V)$ and $n\geq 0$.

\bigskip

\noindent
(ii) Consider the Lie subalgebra $\mathfrak{a}\subseteq\cG$ generated by the elements $x_\ga\ot 1$ and $x_\beta\ot 1$.  Let $v\in V$ be a nonzero vector of weight $\lambda$, and let $W\subseteq V$ be the $\mathfrak{a}$-submodule generated by $v$.  The Lie algebra $\mathfrak{a}$ is clearly finite dimensional, and $W$ is a finitely generated $\mathfrak{a}$-module.  By \cite[Corollary 2.7]{Fe90},
$$\mathfrak{a}[W]=\{x\in\mathfrak{a}\ :\ x\hbox{\ acts locally finitely on\ }W\}$$
is a Lie subalgebra of $\mathfrak{a}$.  Since $\ga +\beta$ is a root, the element $x_{\ga+\beta}\ot 1$ is a nonzero multiple of $[x_\ga\ot 1,x_\beta\ot 1]\in\mathfrak{a}[W].$  Therefore $x_{\ga+\beta}\ot 1$ cannot act injectively on $W$, so it acts nilpotently everywhere on $V$, and $\ga+\beta\in\Phi^f$.

\bigskip

\noindent
(iii) Suppose $\ga\in\Phi^f\cap-\Phi^f$ and $\beta\in\Phi^i\cap-\Phi^i$.  If $\ga+\beta\in\Phi^f$, then by Part (ii),
$$\beta=-\ga+(\ga+\beta)\in\Phi^f,$$
a contradiction.  Likewise, if $\ga+\beta\in\Phi^i$, then by Part (i),
$$\ga=-\beta+(\ga+\beta)\in\Phi^i,$$
another contradiction.  Hence $\ga+\beta\notin\Phi$.\qed

\bigskip

We say that a subset $\mathcal{S}$ of the root system $\Phi$ is {\em closed} if $\ga+\beta\in\mathcal{S}$ whenever $\ga,\beta\in\mathcal{S}$ and $\ga+\beta\in\Phi$.

\begin{lemma}\label{lemma4}
\begin{enumerate}
\item[{\rm (i)}] The sets $P=\Phi^f\cup-\Phi^i$ and $-P=(-\Phi^f)\cup\Phi^i$ are closed.
\item[{\rm (ii)}] The sets $\Phi^i\cap-\Phi^i$, $\Phi^f\cap-\Phi^f$, and $P\cap-P$ are root subsystems of $\Phi$.
\end{enumerate}
 
\end{lemma}

\noindent
\proof (i) Suppose that $\ga,\beta\in P$ and $\ga+\beta\in\Phi$.  If $\ga,\beta\in\Phi^f$ or $\ga,\beta\in-\Phi^i$, then this follows from Lemma \ref{lemma3}(i) or \ref{lemma3}(ii).  Otherwise, we may assume that $\ga\in\Phi^f\setminus-\Phi^i$ and $\beta\in(-\Phi^i)\setminus\Phi^f$.  Then $\ga\in\Phi^f\cap-\Phi^f$ and $\beta\in\Phi^i\cap-\Phi^i$, which is impossible by Lemma \ref{lemma3}(iii), so this case never occurs.

\bigskip

\noindent
(ii) This is an immediate consequence of \cite[Proposition VI.1.23]{bourbaki}.\qed

\begin{corollary}
The subspaces $\mathfrak{l}$, $\mathfrak{n_\pm}$, and $\mathfrak{p}$ are Lie subalgebras of $\fg$. \qed
\end{corollary}


We denote the corresponding current algebras by $\cL=\mathfrak{l}\ot S$, $\mathcal{N}_\pm=\mathfrak{n}_\pm\ot S$, and $\mathcal{P}=\mathfrak{p}\ot S$.


\begin{remark}
Since $\Phi$ is the disjoint union of $\Phi^f$ and $\Phi^i$, the following properties are obvious:
\begin{enumerate}
\item[{\rm (i)}] $P\cup-P=\Phi$
\item[{\rm (ii)}] $P\setminus-P\subseteq\Phi^f$
\item[{\rm (iii)}] $(-P)\setminus P\subseteq\Phi^i$.
\end{enumerate} 
In particular, $P$ is a {\em parabolic} subset of the root system $\Phi$.
\end{remark}
 
\begin{proposition}\label{proposition1}
\begin{enumerate}
\item[{\rm (i)}] The subalgebra $\mathfrak{l}\subseteq\fg$ is reductive, with Cartan subalgebra $\fh$ and root system $P\cap-P$.  More explicitly, $\mathfrak{l}=A(\mathfrak{l})\oplus Z(\mathfrak{l})$, where  $A(\mathfrak{l})=\hbox{Span}\{h_\ga : \ga\in P\cap-P\}\oplus\bigoplus_{\ga\in P\cap-P}\fg_\ga$ is semisimple, $Z(\mathfrak{l})$ is the centre of $\mathfrak{l}$, and $\mathfrak{h}=Z(\mathfrak{l})\oplus\hbox{Span}\{h_\ga : \ga\in P\cap-P\}.$
\item[{\rm (ii)}] Let $W$ be a simple admissible module for $\mathcal{L}=\mathfrak{l}\ot S$, with respect to the Cartan subalgebra $\fh$.  Then
$$W\cong U\ot k_\chi,$$
where $U$ is a simple admissible module for $A(\mathfrak{l})\ot S$, and $k_\chi$ is the $1$-dimensional module for $Z(\mathcal{L})=Z(\mathfrak{l})\ot S$ given by a character $\chi:\ Z(\mathcal{L})\rightarrow k$. 
\end{enumerate}
\end{proposition}

\noindent
\proof Part (i) follows from basic properties of root systems, as seen in \cite[chapitre VI]{bourbaki}, for instance.  Part (ii) can be proved using the argument in \cite[Proposition 3.4]{curr}.\qed

\bigskip

Let $Q$ be the root lattice $\hbox{Span}_\mathbb{Z}\Phi$, and $Q_0$ the sublattice generated by $P\cap-P$.  Recall that $V$ is a simple Harish-Chandra module for $\mathcal{G}$.  If $W\subseteq V$ is a simple $\mathcal{L}$-submodule, then the weights of $W$ differ only by elements of $Q_0$.  Such modules can be interpreted as $\mathcal{P}$-modules by defining the action of the nilradical $\mathcal{N}_+$ to be trivial.

\begin{proposition}\label{proposition2}
Let $W$ be a weight module for $\mathcal{L}$ whose weights differ only by elements of $Q_0$.
\begin{enumerate}
\item[{\rm (i)}] The generalised Verma module $$M_\mathcal{P}(W)=\mathcal{U}(\mathcal{G})\ot_{\mathcal{U}(\mathcal{P})}W$$ has a unique submodule $N_\mathcal{P}(W)$ which is maximal among all submodules having trivial intersection with $W$.
\item[{\rm (ii)}] The quotient $L_\mathcal{P}(W)=M_\mathcal{P}(W)/N_\mathcal{P}(W)$ is a simple $\cG$-module if and only if $W$ is a simple $\mathcal{L}$-module.
\item[{\rm (iii)}] The space $L_\mathcal{P}(W)^\mathcal{N_+}$ of $\mathcal{N}_+$-invariants is precisely $W$.
\end{enumerate}
\end{proposition}

\noindent
\proof 
(i) By the Poincar\'e-Birkhoff-Witt Theorem, $M_\mathcal{P}(W)$ can be identified with 
$$\mathcal{U}(\mathcal{N}_-)\ot_{\mathcal{U}(\mathcal{P})}W=\left(\mathcal{U}(\mathcal{N}_-)\mathcal{N_-}\ot_{\mathcal{U}(\mathcal{P})}W\right)\oplus W$$ 
as vector spaces or $\mathcal{L}$-modules.  Recalling the characterisation of parabolic sets in \cite[Proposition VI.1.7.20]{bourbaki}, there is a base $\Delta\subseteq\Phi$ of simple roots and a subset $T\subseteq \Delta$ such that $P$ is the union of the set $\Phi^+$ of positive roots with respect to $\Delta$, together with the roots in the $\mathbb{N}$-span of $-T$.  In particular, the unique expression of any $\ga\in(-P)\setminus P$ as a (necessarily negative integer) linear combination of simple roots will contain at least one simple root in $\Delta\setminus T$.  As the weights of $W$ differ only by elements of $Q_0$, the supports of the $\mathcal{L}$-modules $\mathcal{U}(\mathcal{N}_-)\mathcal{\mathcal{N}_-}\ot_{\mathcal{U}(\mathcal{P})}W$ and $W$ are disjoint, and any $\mathcal{G}$-submodule $N\subseteq M_\mathcal{P}(W)$ intersecting $W$ trivially is necessarily contained in $\mathcal{U}(\mathcal{N}_-)\mathcal{\mathcal{N}_-}\ot_{\mathcal{U}(\mathcal{P})}W$.

The sum $N_\mathcal{P}(W)$ of all submodules $N\subseteq M_\mathcal{P}(W)$ such that $N\cap W=0$ is thus also contained in $\mathcal{U}(\mathcal{N}_-)\mathcal{\mathcal{N}_-}\ot_{\mathcal{U}(\mathcal{P})}W$, and is therefore the unique $\mathcal{G}$-submodule which is maximal among all submodules having trivial intersection with $W$.

\bigskip

\noindent
(ii) Let $K$ be a $\mathcal{G}$-submodule of $L_\mathcal{P}(W)$.  We identify $W$ with its isomorphic image in $L_\mathcal{P}(W)$.  If $K\cap W=0$, then $K=0$ by the construction of $N_\mathcal{P}(W)$.  Suppose $W$ is a simple $\mathcal{L}$-module and $K\neq 0$.  Then
$$K\supseteq \mathcal{U}(\mathcal{G})(K\cap W)=\mathcal{U}(\mathcal{G})\cU(\cL)(K\cap W)=\cU(\cG)W=L_\cP(W),$$
so $L_\cP(W)$ is simple.

If $W$ is not simple, then there is a nonzero proper $\cL$-submodule $U\subsetneq W$, and $U_\gamma\subsetneq W_\gamma$ for some weight $\gamma$ of $W$.  By the argument of part (i), the submodule $N_\cP(W)\subseteq M_\cP(W)$ is contained in $\cU(\cN_-)\cN_-\ot_{\cU(\cP)}W$, which has support disjoint from $W$.  It follows that the weight space $L_\cP(W)_\lambda=W_\lambda$ for all $\lambda\in\Pi(W)$.  But then we have a proper inclusion of weight spaces $$L_\gamma=U_\gamma\subsetneq W_\gamma=L_\cP(W)_\gamma,$$ where $L$ is the $\cG$-submodule
$$L=\frac{\left(\cU(\cG)\ot_{\cU(\cP)}U\right)+N_\cP(W)}{N_\cP(W)}\subseteq L_\cP(W).$$
Therefore, $L$ is a nonzero proper submodule of $L_\cP(W)$.

\bigskip

\noindent 
(iii) As in the proof of part(ii), we identify $W$ with its isomorphic image in $L_\cP(W)$.  Suppose $v$ is a weight vector in $L_\cP(W)^{\cN_+}\setminus W$.  Then the argument in part (i) shows that $v$ is necessarily in the vector space $\cU(\cN_-)\cN_-W$.

Let $X\subseteq L_\cP(W)$ be the $\cG$-submodule generated by $v$.  Since $v$ is an $\cN_+$-invariant of $L_\cP(W)$, we see that
\begin{align*} 
X&=\cU(\cG)v\\
&=\cU(\cN_-)\cU(\cL)\cU(\cN_+)v\\
&=\cU(\cN_-)\cU(\cL)v\\
&\subseteq \cU(\cN_-)\cU(\cL)\cU(\cN_-)\cN_-W.
\end{align*}
But $[\cL,\cN_-]\subseteq \cN_-$, so it follows that $X\subseteq \cU(\cN_-)\cN_-W$.  As $\cU(\cN_-)\cN_-W$ and $W$ have disjoint supports, $X\cap W=0$.  By the construction of $N_\cP(W)$, this means that $X=0$ in $L_\cP(W)$.  Hence $v=0$ and $L_\cP(W)^{\cN_+}\subseteq W$.  The reverse inclusion is immediate from the definition of $L_\cP(W)$.\qed

\bigskip

We now prove the main theorem of this section.

\begin{theorem}\label{thm1}
Let $V$ be a simple Harish-Chandra module for $\cG$.  Then one of the following (mutually exclusive) cases occurs:
\begin{enumerate}
\item[{\rm (i)}] $\Phi=\Phi^f$ and $V$ is finite dimensional;
\item[{\rm (ii)}] $\Phi=\Phi^i$ and $V$ is infinite dimensional and admissible, with $\dim V_\lambda=\dim V_\mu$ for all $\lambda,\mu\in\Pi(V)$;
\item[{\rm (iii)}] $\Phi$ is neither $\Phi^f$ nor $\Phi^i$.  Then the corresponding parabolic subalgebra $\cP$ is properly contained in $\cG$, $V^{\cN_+}$ is a simple admissible $\cL$-module, and $V\cong L_\cP(V^{\cN_+})$.
\end{enumerate}

\end{theorem}

\noindent
\proof {\bf (i)} Suppose that $\Phi=\Phi^f=\{\ga_1,\ldots,\ga_r\}$.  By the PBW Theorem,
$$V=\cU(\fg_{\ga_r}\ot S)\cU(\fg_{\ga_{r-1}}\ot S)\cdots\cU(\fg_{\ga_1}\ot S)\cU(\fh\ot S)V_\gl,$$
for any $\gl\in\Pi(V)$.  The algebra $\cU(\fh\ot S)$ stabilises the weight space $V_\gl$, and the $\ga_1$-string $\mathcal{S}_1(\gl)=\Pi(V)\cap\{\gl+n\ga_1\ :\ n\geq 0\}$ is finite by Proposition \ref{Phi^f-prop}(ii).  Thus $(\fg_{\ga_1}\ot S)\cU(\fh\ot S)V_\gl$ has only finitely many nonzero $\fh\ot 1$-weight spaces.  Similarly, $(\fg_{\ga_2}\ot S)(\fg_{\ga_1}\ot S)\cU(\fh\ot S)V_\gl$ has finite support, since the $\ga_2$-string $\mathcal{S}_2(\mu)$ is finite for each of the (finitely many) weights $\mu$ of $(\fg_{\ga_1}\ot S)\cU(\fh\ot S)V_\gl$.  By induction on $r$, we see that $V$ has only finitely many nonzero weight spaces, each of which is finite dimensional, since $V$ is a Harish-Chandra module.

\medskip

\noindent
{\bf (ii)} Suppose $\Phi=\Phi^i$, $\gl\in\Pi(V)$, and $\ga\in\Phi$.  The element $x_\ga\ot 1$ acts injectively on $V$, so in particular, the map 
$$x_\ga\ot 1:\ V_\gl\longrightarrow V_{\gl+\ga}$$
is injective and $\dim V_\gl\leq\dim V_{\gl+\ga}$.  Similarly, $-\ga\in\Phi$, so $x_{-\ga}\ot 1:\ V_{\gl+\ga}\longrightarrow V_\gl$
is injective and $\dim V_\gl=\dim V_{\gl+\ga}$.

Since $V$ is simple, $V_\gl$ generates $V$, so $\Pi(V)=\gl+Q$.  The module $V$ is thus infinite dimensional with $\dim V_\mu=\dim V_\gl$ for all $\mu\in\Pi(V)$.

\medskip

\noindent
{\bf (iii)} We now consider the case where $\Phi$ is neither $\Phi^f$ nor $\Phi^i$.  Suppose that $P=\Phi$.  If $\ga\in\Phi^f$, then 
$$-\ga\in\Phi\setminus-\Phi^i=P\setminus-\Phi^i=\left(\Phi^f\cup-\Phi^i\right)\setminus -\Phi^i\subseteq\Phi^f.$$
It follows that 
$$\Phi=\left(\Phi^i\cap-\Phi^i\right)\cup\left(\Phi^f\cap-\Phi^f\right).$$

Suppose that $\ga\in\Phi^i\cap-\Phi^i$ and $\beta\in\Phi^f\cap-\Phi^f$.  By Lemma \ref{lemma3}(iii), neither $\ga+\beta$ nor $\ga-\beta$ is a root.  Therefore, $\ga$ and $\beta$ are orthogonal with respect to Killing form on $\fh^*$, \cite[Lemma 9.4]{humphreys-intro}.  By \cite[Proposition VI.1.7.23]{bourbaki}, $\Phi^i\cap-\Phi^i$ and $\Phi^f\cap-\Phi^f$ are root subsystems of $\Phi$.  But this is an orthogonal decomposition of $\Phi$ into nonempty subsystems, contradicting the simplicity of $\fg$.  Hence $\Phi=\Phi^i\cap-\Phi^i$ or $\Phi=\Phi^f\cap-\Phi^f$, so $\Phi=\Phi^i$ or $\Phi=\Phi^f$, another contradiction.  Hence $P\neq \Phi$, and $\cP$ is a proper subalgebra of $\cG$.

As $[\mathcal{L},\cN_+]\subseteq \cN_+$, it is clear that $V^{\cN_+}$ is an $\cL$-submodule of $V$.  Since $\Phi\neq \Phi^f$, we see that $V\neq 0$.  Let $v_0\in V$ be a nonzero element of weight $\gl$.  By a PBW argument as in Case (i), the $\cN_+$-submodule $\cU(\cN_+)v_0$ is finite dimensional.  Here we use the fact that $\cN_+$ is spanned by elements $x_\ga\ot s$, where $\ga\in P\setminus  -P\subseteq\Phi^f$ and $s\in S$.  By \cite[Proposition VI.1.7.20]{bourbaki}, there is a base $\Delta$ of $\Phi$, with respect to which all roots of $P\setminus-P$ are positive.  There is thus a highest weight in $\cU(\cN_+)v_0$, that is, a weight $\mu\in\Pi(V)$ such that $\left(\cU(\cN_+)v_0\right)\mu\neq 0$, but $\left(\cU(\cN_+)v_0\right)\ga= 0$ for all $\ga>\mu$.  But then $\cN_+$ acts trivially on 
$\left(\cU(\cN_+)v_0\right)_\mu$, so $V^{\cN_+}$ contains $\left(\cU(\cN_+)v_0\right)\mu$ and $V^{\cN_+}\neq 0$.  

The induced module $M_\cP(V^{\cN_+})=\cU(\cN_-)\ot_{\cU(\cP)}V^{\cN_+}$ is free as an $\cN_-$-module, so there is a well-defined $\cG$-module homomorphism 
$$\phi:\ M_\cP(V^{\cN_+})\longrightarrow V$$
given by $\phi(u\ot v)=uv$ for all $u\in\cU(\cN_-)$ and $v\in V^{\cN_+}$.  The map $\phi$ is nonzero since its image contains $V^{\cN_+}$, so is surjective by the simplicity of $V$.  The simple module $V$ is thus the quotient of $M_\cP(V^{\cN_+})$ by $\ker\phi$.

By the proof of Proposition \ref{proposition2}(i), the $\cL$-modules $V^{\cN_+}$ and $\cU(\cN_-)\cN_-\ot_{\cU(\cP)}\left(V^{\cN_+}\right)\subseteq M_\cP(V^{\cN_+})$ have disjoint support.  The map $\phi$ is clearly an $\cL$-module homomorphism and restricts to an injection on $V^{\cN_+}$, while preserving weight spaces.  It follows that $\ker\phi\cap\left( V^{\cN_+}\right)=0$, so $\ker\phi\subseteq N_\cP(V^{\cN_+})$ by Proposition \ref{proposition2}(i).  The submodule $N_\cP(V^{\cN_+})$ is thus a maximal submodule of $M_\cP(V^{\cN_+})$.  Therefore, $L_\cP(V^{\cN_+})\cong V$ is a simple $\cG$-module, and by Proposition \ref{proposition2}(ii), $V^{\cN_+}$ is a simple $\cL$-module.

It remains only to show that $V^{\cN_+}$ is admissible.  By definition, $P\cap-P=\left(\Phi^i\cap-\Phi^i\right)\cup\left(\Phi^f\cap-\Phi^f\right)$, and as we have noted above, this is an orthogonal decomposition of $P\cap-P$ into a disjoint union of two root subsystems of $\Phi$.  In particular, $\mathfrak{l}$ decomposes into a Lie algebra direct sum
$$\mathfrak{l}=\mathfrak{z}\oplus\fg^i\oplus\fg^f,$$
for some abelian Lie algebra $\mathfrak{z}$ and semisimple Lie algebras $\fg^i$ and $\fg^f$, with root systems $\Phi^i\cap-\Phi^i$ and $\Phi^f\cap-\Phi^f$, respectively.

The simple $\cL$-module $V^{\cN_+}$ thus factors into a tensor product:
$$V^{\cN_+}\cong W^0\ot W^i\ot W^f,$$
where $W^0$, $W^i$, and $W^f$ are simple Harish-Chandra modules for the current algebras $\mathfrak{z}\ot S$, $\fg^i\ot S$, and $\fg^f\ot S$, respectively.  The $\mathfrak{z}\ot S$-module $W^0$ is $1$-dimensional, and by Case (i), $W^f$ is also finite dimensional.  If $\fg^i=0$, then $W^i$ can be taken to be the $1$-dimensional trivial module; otherwise, $W^i$ is infinite dimensional and admissible by Case (ii).  Therefore, $V^{\cN_+}\cong W^0\ot W^i \ot W^f$ is an admissible $\cL$-module.\qed

\bigskip

\begin{remark}\label{summary-thm1} Every simple Harish-Chandra module for $\cG$ is thus the irreducible quotient $L_\cP(W)$ of a generalised Verma module $M_\cP(W)$ induced from a simple admissible module $W$ for a Levi subalgebra $\cL\subseteq\cG$.  The converse is non-trivial, as the module $M_\cP(W)$ will clearly have infinite dimensional weight spaces whenever $S$ is infinite dimensional as a vector space.  In the following section, we use exp-polynomial functions and a careful analysis of raising operators to prove that if $W$ is simple and admissible, then $L_\cP(W)$ is a simple Harish-Chandra module for $\cG$.  An alternative proof based on evaluation modules will be presented in Section 4.


\end{remark}


\section{Finite dimensionality of weight spaces}

For any set $T$, let $T^\bullet$ be the set of all finite sequences $(t_1,\ldots,t_p)$, $p\geq 0$, of elements $t_i\in T$.  The {\em length} of an element $\und{t}=(t_1,\ldots,t_p)\in T^\bullet$ is defined to be $p$ and is denoted by $|\und{t}|$.  When $T$ is an additive semigroup, we write $t=t_1+\cdots+t_p$ for the sum of the components of any element $\und{t}=(t_1,\ldots,t_p)\in T^\bullet$; when $T$ is a multiplicative monoid, multi-index notation is used for exponentiation: $\und{t}^{\und{n}}=t_1^{n_1}t_2^{n_2}\cdots t_p^{n_p},$ for each $\und{t}\in T^\bullet$ and $\und{n}\in\mathbb{N}^\bullet$ with $|\und{t}|=|\und{n}|=p$.  Note that we allow $|\und{t}|$ to be zero, in which case $t=0$ and $\und{t}^{\und{n}}=1$, in the additive and multiplicative cases, respectively.  Fix nonzero generators $s_1,\ldots,s_d$ for the finitely generated $k$-algebra $S$.  If $\und{b}=(b_{11},\ldots,b_{1d},b_{21},\ldots,b_{2d},\ldots,b_{\ell 1},\ldots,b_{\ell d})\in\mathbb{N}^{\ell d}$ and $\und{\ga}=(\ga_1,\ldots,\ga_\ell)\in\Phi^\bullet$, then we write 
$$x_{\und{\ga}}(\und{b})=x_{\ga_1}(\und{s}^{\und{b_1}})x_{\ga_2}(\und{s}^{\und{b_2}}) \cdots x_{\ga_\ell}(\und{s}^{\und{b_\ell}})\in\cU(\cG),$$
where  $x_{\ga_i}(\und{s}^{\und{b_i}})=x_{\ga_i}\ot s_1^{b_{i1}}s_2^{b_{i2}}\cdots s_d^{b_{id}}\in\cG$.  To simplify the exposition, we will implicitly assume that $|\und{t}|=|\und{n}|$ and $|\und{b}|=d|\und{\ga}|$ whenever we write $\und{t}^{\und{n}}$ or $x_{\und{\ga}}(\und{b})$, respectively.

 By \cite[Proposition VI.1.7.20]{bourbaki}, the parabolic subsets $P\subseteq\Phi$ are precisely those subsets for which there is a base $\Delta=\{\delta_1,\ldots,\delta_m,\epsilon_1,\ldots,\epsilon_n\}$ of simple roots of $(\fg,\fh,\Phi)$ such that $P=\Phi\cap\left(\hbox{Span}_{\mathbb{N}}\Delta\oplus\hbox{Span}_{\mathbb{N}}\{-\delta_1,\ldots,-\delta_m\}\right)$.  In particular, we can define the {\em (parabolic) height} of an element $\mu$ of the root lattice $Q=\hbox{Span}_{\mathbb{Z}}\Phi$ to be
$$\hgt\mu=\hbox{ht}_P(\mu):=\sum_{j=1}^n d_j,$$
where
$$\mu=\sum_i c_i\delta_i+\sum_j d_j\epsilon_j$$
is the (unique) expression of $\mu$ in terms of the base $\Delta$.  For any $\und{\mu}\in\Phi^\bullet$, we define $\hgt\und{\mu}:=\hbox{ht}_P(\und{\mu})$ to be the sum $\sum_i\hgt\mu_i$ of the heights of its components $\mu_i$.  Note that if $P$ is the set of positive roots with respect to $\Delta$, this definition reduces to the ordinary notion of height.

Let $P$ be a parabolic subset of $\Phi$, and let $\mathfrak{l}=\fh\oplus\sum_{\ga\in P\cap-P}\fg_\ga$ and $\mathfrak{p}=\fh\oplus\sum_{\ga\in P}\fg_\ga$ be the corresponding Levi and parabolic subalgebras of $\fg$.  As in the previous section, we write $\cG=\fg\ot S$, $\cL=\mathfrak{l}\ot S$, and $\cP=\mathfrak{p}\ot S$ for the associated current algebras.  Let $W$ be a simple admissible $\cL$-module.  We will show that $L_\cP(W)=M_\cP(W)/N_\cP(W)$ is a (simple) Harish-Chandra module for $\cG$.

For convenience, we extend the definition of height to the support $\Pi(M)$ of the module $M=M_\cP(W)$.  Fix a weight $\gl\in\Pi(W)$.  For any $\gamma\in\Pi(M)$, note that $\gamma-\gl$ is in the root lattice $Q$, so we can define $\hgt_M\gamma$ to be $\hgt\!(\gamma-\gl)$.  As $\gl_1-\gl_2\in\hbox{Span}_\mathbb{Z}(P\cap-P)$ for all $\gl_1,\gl_2\in\Pi(W)$, we see that $\hgt\!(\gl_1-\gl_2)=0$, and this definition is independent of the choice of $\gl\in\Pi(W)$.  Similarly, we define the {\em height} $\hgt v$ of a weight vector $v\in M$ to be the height $\hgt_M \mu$ of its weight $\mu$.

\begin{lemma}\label{lemma1-fin}
Let $v\in M_\cP(W)$ be a nonzero vector of weight $\gamma$.  Then $v\in N_\cP(W)$ if $x_{\und{\beta}}(\und{b})v=0$ whenever $\und{\beta}\in(P\setminus-P)^\bullet$, $\und{b}\in\mathbb{N}^\bullet$, and {\em $\hgt\und{\beta}=-\hgt_M\gamma$}.  Moreover, the number of $\und{\beta}\in(P\setminus-P)^\bullet$ with {\em $\hgt\und{\beta}=-\hgt_M\gamma$} is finite.
\end{lemma}

\noindent
\proof For $w\in W$, the height of any element $x_{\und{\ga}}(\und{a})w$ is simply $\hgt\und{\ga}$, since $\hgt w=0$.  Since $M_\cP(W)=\cU(\cN_-)\ot_{\cU(\cP)}W$ and $\hgt\mu<0$ for all $\mu\in(-P)\setminus P$, it follows that every element of $M_\cP(W)$ has non-positive height.  In particular, $x_{\und{\beta}}(\und{b})v=0$ whenever $\und{\beta}\in(P\setminus-P)^\bullet$, $\und{b}\in\mathbb{N}^\bullet$, and $\hgt\und{\beta}>-\hgt_M\gamma$.

As every element of $(-P)\setminus P$ has negative height and every element of $P\setminus-P$ has positive height, the only vectors of height $0$ in $\cU(\cG)v=\cU(\cN_-)\cU(\cL)\cU(\cN_+)v$ are thus $\cU(\cL)$-linear combinations of those of the form $x_{\und{\beta}}(\und{b})v$, for some $\und{\beta}\in(P\setminus-P)^\bullet$ and $\und{b}\in\mathbb{N}^\bullet$ with $\hgt\und{\beta}=-\hgt_M\gamma$.  Since every element of $W$ has height $0$, we see that $\cU(\cG)v\cap W=0$ if $x_{\und{\beta}}(\und{b})v=0$ for all $\und{\beta}\in(P\setminus-P)^\bullet$ and $\und{b}\in\mathbb{N}^\bullet$ with $\hgt\und{\beta}=-\hgt_M\gamma$.  In particular, the $\cG$-submodule $\cU(\cG)v$ is contained in $N_\cP(W)$.

Since every element $\mu\in P\setminus-P$ has strictly positive height, the length of each $\und{\beta}\in (P\setminus-P)^\bullet$ for which $\hgt\und{\beta}=-\hgt_M\gamma$ is bounded by $-\hgt_M\gamma$, so the number of $\und{\beta}\in (P\setminus-P)^\bullet$ with $\hgt\und{\beta}=-\hgt_M\gamma$ is finite.\qed

\bigskip

The algebra of {\em exp-polynomial functions}\footnote{The original definition of exp-polynomial was given by Billig and Zhao \cite{BiZh04} for $\mathbb{Z}^n$-graded Lie algebras, generalising previous work of Berman and Billig on toroidal Lie algebras \cite{BeBi99}.  Our definition differs somewhat from theirs in order to accommodate our non-graded context.} in $r$ variables $a_1,\ldots,a_r$ is the algebra of functions $f:\ \mathbb{N}^r\rightarrow k$, generated by the {\em exponential functions} $\und{a}\mapsto\und{b}^{\und{a}}$ and the {\em polynomial functions} $\und{a}\mapsto\und{a}^{\und{c}}$ for $\und{b}\in k^r$ and $\und{c}\in\mathbb{N}^r$.  By convention, $0^0$ is defined to be $1$.  Any exp-polynomial function in $m+n$ variables $a_1,\ldots,a_m,b_1,\ldots,b_n$ can be expanded as a finite sum
$$f(\und{a},\und{b})=\sum_p f_p(\und{a})\und{c_p}^{\und{b}}\und{b}^{\und{d_p}},$$
for some $\und{c_p}\in k^n$, $\und{d_p}\in\mathbb{N}^n$, and exp-polynomial functions $f_p(\und{a})$ in $m$ variables $a_1,\ldots,a_m$.

Let $\{x_1,\ldots,x_m\}$ be a basis for the Lie algebra $\fg$, which is homogeneous with respect to the root grading, and has structure constants $c_{ij}^\ell\in k$ given by $[x_i,x_j]=\sum_{\ell=1}^mc_{ij}^\ell x_\ell.$
In the commutation relations
\begin{equation}\label{eqn1}
[x_i(\und{a}),x_j(\und{b})]=\sum_{\ell=1}^mc_{ij}^\ell x_\ell(\und{a}+\und{b})
\end{equation}
governing the multiplication of $\cU(\cG)$, the constant coefficients $c_{ij}^\ell$ are obviously exp-polynomial functions of the $2d$ variables $(\und{a},\und{b})=(a_1,\ldots,a_d,b_1,\ldots,b_d)$, where $x_i(\und{a})$ denotes $x_i\ot \und{s}^{\und{a}}$, for instance.

By \cite{curr}, the simple admissible $\cL$-modules $W$ can be written as tensor products
$$W\cong W_1(M_1)\ot\cdots\ot W_q(M_q)$$
of evaluation modules $W_1(M_1),\ldots,W_q(M_q)$, where $W_1,\ldots, W_q$ are simple $\fl$-modules, $M_1,\ldots, M_q$ are distinct maximal ideals of $S$, and the tensor product $W_1\ot\cdots\ot W_q$ is admissible.  That is, $W$ is isomorphic to the $\cL$-module $W_1\ot\cdots\ot W_q$, with action given by
$$(x\ot s).(w_1\ot\cdots\ot w_q)=\sum_{i=1}^q s(M_i)w_1\ot\cdots\ot xw_i\ot\cdots\ot w_q$$
for each $x\in\fl$, $s\in S$, and $w_i\in W_i$.  The element $s(M_i)$ of the algebraically closed field $k$ is the reduction of $s$ modulo $M_i$:
$$s(M_i)+M_i=s+M_i\in S/M_i\cong k.$$
Reduction is a ring homomorphism, so writing $b_{ij}=s_j(M_i)$ for each of the generators $s_1,\ldots,s_d$ of $S$, we see that
$$x_\ell(\und{a}).(u_0\ot\cdots\ot u_q)=\sum_{i=1}^q\und{b_i}^{\und{a}}u_0\ot\cdots\ot x_\ell u_i\ot\cdots\ot u_q,$$
for $\und{b_i}=(b_{i1},\ldots,b_{id})$.

Such an action is necessarily {\em exp-polynomial}.  That is, there is a basis $\{w_j\ :\ j\in J\}$ of weight vectors of $W$ and exp-polynomial functions $f_{ij}^\ell$ in $d$ variables $a_1,\ldots,a_d\in\mathbb{N}$, such that
\begin{equation}\label{eqn2-finite}
x_i(\und{a}).w_j=\sum_{\ell\in J} f_{ij}^\ell(\und{a})w_\ell.
\end{equation}
There are only finitely many nonzero functions $f_{ij}^\ell$ for each pair $(i,j)$.  Indeed $x_i$ is a root vector and $w_j$ is a weight vector, so the number of nonzero functions $f_{ij}^\ell$ occuring in the sum is bounded by the maximum dimension of the weight spaces of the admissible $\cL$-module $W$.

\begin{theorem}\label{thm2}
Let $W$ be a simple admissible $\cL$-module.  Then $L_\cP(W)$ is a simple Harish-Chandra module.
\end{theorem}

\noindent
\proof To simplify notation in the proof, we write $M=M_\cP(W)$ and $L=L_\cP(W)$.  Fix a basis $\{w_j\ :\ j\in J\}$ of weight vectors for $W$, as in the previous paragraph, and let $\gamma\in\Pi(M)$.  To avoid ambiguity, we write $\overline{m}$ for the image of each $m\in M$ under the quotient map $M\rightarrow L$.

By Proposition \ref{proposition2}(ii), we need only show that the weight space $L_\gamma$ is finite dimensional.  This weight space is spanned by the set of elements $\overline{x_{\und{\ga}}(\und{a})w_j}$ satisfying the following conditions:
\begin{enumerate}
\item[{\rm (C1)}] $\und{\ga}\in((-P)\setminus P)^\bullet$,
\item[{\rm (C2)}] $\und{a}=(\und{a_1},\ldots,\und{a_\ell})=(a_{11},\ldots,a_{1d},\ldots,a_{\ell 1},\ldots,a_{\ell d})\in\mathbb{N}^{\ell d}$, where $\ell=|\und{\ga}|$,
\item[{\rm (C3)}] $\gamma-\ga$ is the weight of $w_j$.
\end{enumerate}
In particular, $\hgt \und{\ga}=\hgt_M\gamma$, so there are only finitely many such $\und{\ga}$, since each element of $(-P)\setminus P$ has strictly negative height.  This also means that there are only finitely many such $w_j$ occuring, since the weight of each allowable $w_j$ is $\gamma-\ga$ for one of the finitely many $\und{\ga}$, and the dimension of the weight spaces $W_{\gamma-\ga}$ is always finite.


By Lemma \ref{lemma1-fin}, a linear combination $v=\sum_{\und{\ga},\und{a},j} c_{\und{\ga},\und{a},j} \overline{x_{\und{\ga}}(\und{a})w_j}$ of such elements is zero if 
\begin{equation}\label{etoile}
\sum_{\und{\ga},\und{a},j} c_{\und{\ga},\und{a},j}x_{\und{\beta}}(\und{b}) x_{\und{\ga}}(\und{a})w_j=x_{\und{\beta}}(\und{b})v=0
\end{equation}
 for all $\und{\beta}\in(P\setminus-P)^\bullet$ and $\und{b}\in\mathbb{N}^\bullet$ such that $\hgt\und{\beta}=-\hgt_M\gamma$.  The number of possible values of $\und{\beta},\und{\ga},$ and $j$ in (\ref{etoile}) is thus finite and bounded by a number that depends only on $\gamma$.  

Combining (\ref{eqn1}) and (\ref{eqn2-finite}), we see that there exist exp-polynomial functions $f_{\und{\beta},\und{\ga},j}^\ell(\und{b},\und{a})$ in $d|\und{\ga}|+d|\und{\beta}|$ variables $b_{rs},a_{ts}$ with $1\leq r\leq|\und{\beta}|$, $1\leq s\leq d$, and $1\leq t\leq|\und{\ga}|$, such that
$$x_{\und{\beta}}(\und{b}) x_{\und{\ga}}(\und{a})w_j=\sum_{\ell}f_{\und{\beta},\und{\ga},j}^\ell(\und{b},\und{a}) w_\ell.$$
Equation (\ref{etoile}) now becomes 
$$\sum_{\und{\ga},\und{a},j} c_{\und{\ga},\und{a},j}\sum_{\ell}f_{\und{\beta},\und{\ga},j}^\ell(\und{b},\und{a}) w_\ell=0.$$
Since the number of $\und{\beta},\und{\ga},$ and $j$ that can occur is finite and bounded by a constant that depends only on $\gamma$, the number of $\ell$ that can occur is also finite and bounded by a constant that depends only on $\gamma$.  However, the number of possible $\und{b},\und{a}\in\mathbb{N}^\bullet$ is infinite.

Expanding each of the (finitely many) $f_{\und{\beta},\und{\ga},j}^\ell(\und{b},\und{a})$ in the variables $\und{b}$, we obtain
\begin{equation}\label{csillagok}
f_{\und{\beta},\und{\ga},j}^\ell(\und{b},\und{a})=\sum_p h_p(\und{a})\und{c_p}^{\und{b}}\und{b}^{\und{d_p}},
\end{equation}
for some $\und{c_p}\in k^{md}$, $\und{d_p}\in\mathbb{N}^{md}$, and finitely many exp-polynomial functions $h_p$.  The particular functions $h_p$ that occur in (\ref{csillagok}) depend on $\und{\beta},\und{\ga},j,$ and $\ell$, but the total number of such functions is bounded by a number $A$ that depends only on $\gamma$.  

To show that $L_\gamma$ is finite dimensional, it thus suffices to show that for any sufficiently large (finite) set $\mathcal{I}$ of triples $(\und{\ga},\und{a},j)$ satisfying (C1)-(C3), there are always nontrivial coefficients $c_{\und{a}}\in k$ such that
\begin{equation}\label{3etoiles}
\sum_{\und{a}}c_{\und{a}}h_p(\und{a})=0
\end{equation}
for each of the (at most) $A$ possible values of the index $p$.  Since the number of allowable values of $\und{\ga}$ and $j$ is finite and bounded, we have only to choose $\mathcal{I}$ of sufficiently large cardinality to guarantee that the number $B$ of values of $\und{a}$ occuring in triples $(\und{\ga},\und{a},j)\in\mathcal{I}$ is strictly greater than $A$.  Then the linear system (\ref{3etoiles}) consists of at most $A$ linear equations in $B>A$ variables $c_{\und{a}}$.  There is thus a nontrivial solution, and the set of vectors $\left\{\overline{x_{\und{\ga}}(\und{a})w_j}\ :\ (\und{\ga},\und{a},j)\in\mathcal{I}\right\}$
is linearly dependent whenever $\mathcal{I}$ is sufficiently large.  In particular, the weight space $L_\gamma$ is finite dimensional.\qed

\bigskip

Theorems \ref{thm1} and \ref{thm2} can be summarized as follows:

\begin{corollary}\label{thm3} The simple Harish-Chandra modules for the current algebra $\cG$ are precisely the modules $L_\cP(W)$ obtained via parabolic induction from simple admissible modules $W$ of Levi subalgebras.\end{corollary}\qed



\section{Evaluation modules}

Let $\cL=\fl\ot S$, where $\fl$ is the generalised Levi subalgebra $\displaystyle{\fh\ \oplus\sum_{\ga\in P\cap -P}\fg_\ga}$ defined by a parabolic subset $P\subseteq\Phi$.  By \cite{curr}, the simple admissible $\cL$-modules $W$ are isomorphic to tensor products
$$W\cong W_1(M_1)\ot\cdots\ot W_q(M_q)$$
of evaluation modules $W_1(M_1),\ldots,W_q(M_q)$, where $W_1,\ldots, W_q$ are simple admissible $\fl$-modules, $M_1,\ldots, M_q$ are distinct maximal ideals of $S$, and the tensor product $W_1\ot\cdots\ot W_q$ is admissible.


The generalised Verma module $M_\cP(W)$ is free as a $\cU(\cN_-)$-module, where the subalgebra $\cN_-=\mathfrak{n}_-\ot S$ is defined as in Section 2.  Let
$$\Delta:\ \cU(\cN_-)\longrightarrow\cU(\cN_-)^{\ot q}$$
be the unital associative algebra homomorphism defined on elements of $\cN_-$ by
$$\Delta(x)=\sum_{i=1}^qx_i,$$
for simple tensors $x_i=1\ot \cdots\ot x\ot\cdots \ot 1$ with $x$ in the $i$th position and $1$ elsewhere.  This generalised coproduct $\Delta$ induces a well-defined $\cG$-module homomorphism
\begin{eqnarray*}
\delta:\ M_\cP(W)&\longrightarrow&\bigotimes_{i=1}^qM_\cP(W_i(M_i))\\
u\ot_{\cU(\cP)}(w_1\ot\cdots\ot w_q)&\longmapsto&\Delta(u)(\ol{w_1}\ot\cdots\ot\ol{w_q})
\end{eqnarray*}
where $u\in\cU(\cN_-)$, $\ol{w_i}=1\ot_{\cU(\cP)}w_i\in M_\cP(W_i(M_i))$, and $(u_1\ot\cdots\ot u_q)(\ol{w_1}\ot\cdots\ot\ol{w_q})$ is defined as $u_1\ol{w_1}\ot\cdots\ot u_q\ol{w_q}$ for each $u_i\in\mathcal{U}(\mathcal{N}_-)$ and $w_i\in W_i(M_i)=W_i$.

As $W_i(M_i)$ is an evaluation representation at $M_i\in\Max\,S$, there is also a $\cG$-module evaluation homomorphism
\begin{eqnarray*}
ev_i:\ M_\cP(W_i(M_i))&\longrightarrow&M_\fp(W_i)(M_i)\\
u\ot_{\cU(\cP)}w_i\,\ &\longmapsto&u(M_i)\ot_{\cU(\fp)}w_i,
\end{eqnarray*}
where $M_\fp(W_i)=\cU(\fg)\ot_{\cU(\fp)}W_i$ is the generalised Verma module for the finite dimensional Lie algebra $\fg$, and $u(M_i)$ is the image of $u\in\cU(\cG)$ under the unital associative algebra homomorphism $\cU(\cG)\rightarrow \cU(\fg)$ extending the evaluation map $x\ot s\mapsto s(M_i)x$ for all $x\ot s\in\cG$.  %
We write $M_\fp(W_i)(M_i)$ for the corresponding evaluation module for the current algebra $\cG$:
\begin{eqnarray*}
\cG\ \ &\longrightarrow&\ \ \ \fg\ \ \longrightarrow\ \ \hbox{End}(W_i)\\
x\ot s&\longmapsto&s(M_i)x.
\end{eqnarray*}
There is then a $\cG$-module homomorphism
\begin{eqnarray*}
\ev:\ \bigotimes_{i=1}^qM_\cP(W_i(M_i))&\longrightarrow&\ \ \bigotimes_{i=1}^qM_\fp(W_i)(M_i)\\
v_1\ot\cdots\ot v_q\ \ &\longmapsto&\ev_1(v_1)\ot\cdots\ot\ev_q(v_q)
\end{eqnarray*}
for all $v_i\in M_\cP(W_i(M_i))$.

The quotient maps $\pi_i:\ M_\fp(W_i)(M_i)\longrightarrow L_\fp(W_i)(M_i)$
from the $M_\fp(W_i)$ to their unique simple quotients $L_\fp(W_i)$ induce a map
\begin{eqnarray*}
\pi:\ \bigotimes_{i=1}^qM_\fp(W_i)(M_i)&\longrightarrow&\ \ \bigotimes_{i=1}^qL_\fp(W_i)(M_i)\\
m_1\ot\cdots\ot m_q\ \ &\longmapsto&\pi_1(m_1)\ot\cdots\ot\pi_q(m_q),
\end{eqnarray*}
and the composition
$$\pi\circ\ev\circ\delta:\ M_\cP(W)\longrightarrow\bigotimes_{i=1}^qL_\fp(W_i)(M_i)$$
is a nonzero $\cG$-module homomorphism to the simple $\cG$-module $\bigotimes_{i=1}^qL_\fp(W_i)(M_i)$.  By uniqueness of the simple quotient  (Proposition \ref{proposition2}(i) and (ii)), $L_\cP(W)\cong\bigotimes_{i=1}^qL_\fp(W_i)(M_i),$ completing the proof of the following proposition:
\begin{proposition}\label{star}
Let $V$ be a simple Harish-Chandra module for $\cG$.  Then there is a parabolic subset $P\subseteq\Phi$, a finite collection $W_1,\ldots,W_q$ of simple admissible modules for $\fl=\fh\oplus\sum_{\ga\in P\cap-P}\fg_\ga$ (with respect to $\fh$), and distinct maximal ideals $M_1,\ldots,M_q\in\Max\,S$, such that $W_1\ot\cdots\ot W_q$ is an admissible $\fl$-module and $V\cong\bigotimes_{i=1}^qL_\fp(W_i)(M_i).\hfill\hbox{\qed}$
\end{proposition}

Conversely, it is clear from Theorem \ref{thm1}, Theorem \ref{thm2}, and Proposition \ref{star} that every module of this form is a simple Harish-Chandra module.  Alternatively, the proposition that follows may be combined with Theorem \ref{thm1} and Proposition \ref{star} to give a second independent proof of Theorem \ref{thm2}.
\begin{proposition}\label{starstar}
Let $P\subseteq\Phi$ be a parabolic subset, and let $W_1,\ldots,W_q$ be simple Harish-Chandra modules for the corresponding Levi subalgebra $\fl$, such that $W_1\ot\cdots\ot W_q$ is a Harish-Chandra module for $\fl$.  Then $L_\fp(W_1)(M_1)\ot\cdots\ot L_\fp(W_q)(M_q)$ is a simple Harish-Chandra $\cG$-module for any collection of distinct maximal ideals $M_1,\ldots,M_q$ of $S$.\hfill\hbox{\qed}
\end{proposition}

\noindent
\proof
As $M_1,\ldots,M_q$ are distinct, the evaluation map
\begin{eqnarray*}
\cG=\fg\ot S&\longrightarrow&\fg\ot(S/M_1\oplus\cdots\oplus S/M_q)=\fg\oplus\cdots\oplus\fg\\
x\ot s&\longmapsto&s(M_1)x\oplus\cdots\oplus s(M_q)x
\end{eqnarray*}
is surjective by the Chinese remainder theorem, and
$L_\fp(W_1)(M_1)\ot\cdots\ot L_\fp(W_q)(M_q)$ is clearly simple.  

Such modules will be Harish-Chandra modules if and only if the corresponding $\fg$-module
$L_\fp(W_1)\ot\cdots\ot L_\fp(W_q)$ has only finite dimensional weight spaces.  But since $W_1,\ldots,W_q$ are simple Harish-Chandra modules for $\fl$, the set of weights $P(W)\subseteq\fh^*$ of $W=W_1\ot\cdots\ot W_q$ lies in a single coset of $Q_0=\hbox{Span}_\Z(P\cap-P).$  The set of weights of $\cU(\mathfrak{n}_-)$ is precisely $C=\hbox{Span}_{\mathbb{N}}((-P)\setminus P)$, and $C\cap Q_0=0$.  This means that for each $\gamma\in\fh^*$, there are at most finitely many pairs $(c,\mu)\in C\times P(W)$ for which $c+\mu=\gamma$.  But $W_1\ot\cdots\ot W_q$ is also a Harish-Chandra module for $\fl$, so for each $\mu\in P(W)$, there are also only finitely many choices of weights $\mu_i$ of $W_i$ for which $\mu_1+\cdots+\mu_q=\mu$.  It follows that $M_\fp(W_1)\ot\cdots\ot M_\fp(W_q)$, and thus  $L_\fp(W_1)\ot\cdots\ot L_\fp(W_q)$, are Harish-Chandra modules.\qed

\bigskip

The main results of this article may be summarized as follows:
\begin{theorem}\label{mainresult}
\begin{enumerate}
\item[{\rm (i)}] The simple Harish-Chandra modules of $\cG=\fg\ot S$ are precisely the modules $L_\cP(W)$ obtained from parabolic induction from simple admissible modules $W$ of generalised Levi subalgebras $\cL=\fl\ot S$, where $\fl=\fh\oplus\sum_{\ga\in P\cap-P}\fg_\ga$ with respect to a parabolic subset $P\subseteq \Phi$ defining the algebras $\fp=\sum_{\ga\in P}\fg_\ga$ and $\cP=\fp\ot S$.
\item[{\rm (ii)}] The simple Harish-Chandra modules $L_\cP(W)$ for $\cG$ are isomorphic to tensor products $L_\fp(W_1)(M_1)\ot\cdots\ot L_\fp(W_q)(M_q)$ of evaluation modules, where $M_1,\ldots,M_q$ are distinct maximal ideals of $S$, $W_1,\ldots,W_q$ are simple admissible $\fl$-modules, and $W_1\ot\cdots\ot W_q$ is an admissible module for $\fl$.  Conversely, every module of this form is a simple Harish-Chandra module.
\item[{\rm (iii)}] Two such modules $\bigotimes_{i=1}^q L_\fp(W_i)(M_i)$ and  $\bigotimes_{j=1}^{q'} L_{\fp'}(W'_j)(M'_j)$ are isomorphic if and only if $q=q'$ and up to a permutation of the tensor factors, $M_i=M_i'$ and $L_\fp(W_i)\cong L_{\fp'}(W'_i)$ for all $i$. 
\end{enumerate}
\end{theorem}

\noindent
\proof Parts (i) and (ii) are Theorem \ref{thm1}, Theorem \ref{thm2}, Proposition \ref{star}, and Proposition \ref{starstar}.  Part (iii) follows easily from the proof of the analogous theorem for the admissible case \cite[Theorem 3.6]{curr}.\qed

\begin{remark}  The classification of simple admissible $\fl$-modules is given in \cite{Ma00}.
\end{remark}

\begin{remark} After a routine reduction to the case where $W_i$ is cuspidal, there is a well-known isomorphism criterion for the $L_\fp(W_i)$ using an action of the small Weyl group $\mathcal{W}(W_i)$.  Namely, $L_\fp(W_i)\cong L_{\fp'}(W'_i)$ if and only if $\fp'=\fp^w$ and $W'_i\cong W_i^w$ for some $w\in\mathcal{W}(W_i)$.  See \cite[Section 1]{Ma00} or \cite[Section 5]{DiMaPe} for notation and details.
\end{remark}

\begin{remark}
By an easy trace argument \cite[Theorem 2.2]{curr}, the central elements of the universal central extension $\widehat{\cG}$ of $\cG$ must act trivially on any simple Harish-Chandra module for $\widehat{\cG}$, so Theorem \ref{mainresult} also describes the simple Harish-Chandra modules of the universal central extension.

An arbitrary central extension of $\cG$ is a Lie algebra direct sum $\widetilde{\cG}\oplus Z$, where $Z$ is an abelian Lie algebra and $\widetilde{\cG}$ is a central quotient of the universal central extension $\widehat{\cG}$.  It follows that the central elements of $\widetilde{\cG}$ act trivially on any simple Harish-Chandra module for $\widetilde{\cG}\oplus Z$, and the elements of $Z$ can act by any central character $\chi:\ Z\rightarrow k$. 
\end{remark}


\end{document}